\documentclass[12pt, a4paper]{article}

\usepackage[utf8]{inputenc}
\usepackage[T1]{fontenc}
\usepackage[english]{babel}
\usepackage{indentfirst}
\usepackage{amsmath}
\usepackage{amsfonts}
\usepackage{amssymb}
\usepackage{amsthm}
\usepackage{breqn}
\usepackage{graphicx}
\usepackage{babelbib}
\usepackage{lmodern} \normalfont
\usepackage{color}
\usepackage{dsfont}

\newtheorem{thm}{Theorem}[section]
\newtheorem{prop}[thm]{Proposition}
\newtheorem{lem}[thm]{Lemma}
\newtheorem{cor}[thm]{Corollary}

\newcommand{\sV}{{\textsf{V}}}
\newcommand{\sH}{{\textsf{H}}}
\newcommand{\sF}{{\textsf{F}}}
\newcommand{\sK}{{\textsf{K}}}
\newcommand{\sA}{{\textsf{A}}}
\newcommand{\sC}{{\textsf{C}}}
\newcommand{\sD}{{\textsf{D}}}

\newcommand{\IN}{\mathbb{N}}
\newcommand{\IP}{\mathbb{P}}
\newcommand{\IE}{\mathbb{E}}
\newcommand{\Z}{\mathbb{Z}}

\newcommand{\R}{\mathbb{R}}
\newcommand{\F}{{\mathcal{F}}}
\newcommand{\A}{{\mathcal{A}}}
\newcommand{\B}{{\mathcal{B}}}

\newcommand{\I}{{\mathcal{I}}}
\newcommand{\II}{{\mathfrak{I}}}
\newcommand{\M}{{\mathcal{M}}}
\newcommand{\MM}{{\mathfrak{M}}}
\newcommand{\N}{{\mathcal{N}}}

\newcommand{\mS}{{\mathcal{S}}}

\newcommand{\W}{{\mathcal{W}}}
\newcommand{\Y}{{\mathcal{Y}}}

\newcommand{\E}{{\mathcal{E}}}
\newcommand{\capa}{\mathop{\mathrm{cap}}}
\newcommand{\dist}{\mathop{\mathrm{dist}}}
\newcommand{\diam}{\mathop{\mathrm{diam}}}

\newcommand{\ome}{\boldsymbol{\omega}}
\newcommand{\bdelta}{\boldsymbol{\delta}}
\newcommand{\ind}{{\mathds{1}}}
\newcommand{\eqlaw}{\stackrel{\text{\tiny law}}{=}}

\newcommand{\dtv}{\mathop{\mathrm{d}_{\mathrm{TV}}}}
\newcommand{\BB}{{\mathsf B}}

\title{An improved decoupling inequality for random interlacements} 

\date{}

\author{Diego F.\ de Bernardini\thanks{Department of Statistics, Institute of Mathematics, Statistics and Scientific Computation, University of Campinas - UNICAMP, Brazil, e-mails: $\{$bernardini, gallesco, popov$\}$@ime.unicamp.br}
\and Christophe Gallesco\footnotemark[1]
\and Serguei Popov\footnotemark[1]}

\begin{document}

\maketitle

\begin{abstract}
\footnotesize 
In this paper we obtain a decoupling feature of the random interlacements process $\I^u\subset \Z^d$, at level $u$, $d\geq 3$. 
More precisely, we show that the trace of the random interlacements process on two disjoint finite sets, $\sF$ and its translated $\sF+x$, can be coupled with high probability of success, when $\|x\|$ is large, with the trace of a process of independent excursions, which we call the \textit{noodle soup} process.
As a consequence, we obtain an upper bound on the covariance between 
two~$[0,1]$-valued functions depending on the configuration of the random interlacements on $\sF$ and $\sF+x$, respectively. This improves a previous bound obtained 
by Sznitman in~\cite{Szn10}.

\vspace{0.3cm}


\end{abstract}

\section{Introduction and results}\label{intro}

Let $d\geq 3$ and~$\sK_1\subset\Z^d$ be a finite set. Without loss of generality, we assume that $0\in\sK_1$. We consider $\sK_2=\sK_1+\hat{x}$, for $\hat{x}\in\Z^d$ such that $\sK_1\cap\sK_2=\emptyset$, and we denote by $\dist(\sK_1, \sK_2)$ the distance between $\sK_1$ and $\sK_2$,
\begin{align*}
\dist(\sK_1, \sK_2) = \min\{\|x-y\|: x\in\sK_1, y\in\sK_2\},
\end{align*}
where $\|\cdot\|$ is the Euclidean norm. Also, we denote by $\diam(\sK_1)$ the diameter of $\sK_1$
\begin{align*}
\diam(\sK_1) = \max\{\|x-y\|: x,y\in\sK_1\}.
\end{align*}

For a positive real number~$\rho$, we denote by~$\BB_{\rho}(y)$ the 
open discrete ball in~$\Z^d$ centered at~$y\in\Z^d$, of radius~$\rho$, 
with respect to the Euclidean norm, i.e.,
\begin{align*}
\BB_{\rho}(y) = \{x\in\Z^d : \|x-y\| < \rho\}.
\end{align*}
For $R = \frac{\|\hat{x}\|-1}{2}$, we consider $\BB_{R}(0)$. We assume that $\hat{x}$ is large enough such that $\BB_{R}(0)$ and $\BB_{R}(0)+\hat{x}$ are disjoint and $\sK_1\subset \BB_R(0)$. For the sake of brevity, we denote
\begin{align*}
\BB_{R}^1:=\BB_{R}(0) ~\text{ and }~ \BB_{R}^2:=\BB_{R}(0) + \hat{x}.
\end{align*}

In this paper, we will show that the trace of the random interlacements (RI) process at level $u$ on the set $\sK_1\cup \sK_2$, denoted by $\I^u_{\sK_1\cup \sK_2}$, can be coupled with high probability of success (when $\|\hat{x}\|$ is large) with the trace on~$\sK_1\cup \sK_2$, denoted by $\M^u_{\sK_1\cup \sK_2}$, of a process of independent excursions called the \textit{noodle soup} (NS) process, which can be described as follows. Let~$\lambda$ be the expected number of excursions performed by the trajectories of the random interlacements, at level~$u$, between the boundary of the set $\sK_1\cup \sK_2$ and the (external) boundary of the set $\BB_R^1\cup \BB_R^2$. The noodle soup process can be described as a Poisson($\lambda$) number of independent excursions of simple random walk (SRW) between the boundaries of the sets $\sK_1\cup \sK_2$ and $\BB_R^1\cup \BB_R^2$, with the SRW excursions starting at points chosen accordingly to the harmonic measure on $\sK_1\cup \sK_2$. 
This will allow us to obtain an upper bound on the total variation distance between the two processes~$\I^u_{\sK_1\cup \sK_2}$ and~$\M^u_{\sK_1\cup \sK_2}$.

We recall that the total variation distance between two probability measures~$P$ and~$\tilde{P}$ defined on the same $\sigma$-field $\F$ is defined by
\begin{align*}
\dtv(P,\tilde{P}):=\sup_{A\in \F}|P[A]-\tilde{P}[A]|. 
\end{align*}
When dealing with random elements $X$ and $Y$, we will write 
(with a slight abuse of notation) $\text{d}_{\text{TV}}(X,Y)$ 
to denote the total variation distance between the laws of~$X$
and~$Y$.

We can now state the following theorem, which is the main result 
of our paper.
{\thm \label{Maintheo} There exist positive constants $C_1$ and $C_2$, 
depending only on the dimension $d$, such that, if $\dist(\emph{\sK}_1,\emph{\sK}_2) \geq C_1 \diam(\emph{\sK}_1)$, then
\begin{align*}
\dtv(\I^u_{\emph{\sK}_1\cup \emph{\sK}_2},\M^u_{\emph{\sK}_1\cup \emph{\sK}_2}) 
&\leq C_2 \sqrt{u} \frac{\capa(\emph{\sK}_1)^{\frac{3}{2}}}{\dist(\emph{\sK}_1,\emph{\sK}_2)^{d-2}}.
\end{align*}
}

The above result allows to quantify the asymptotic dependence between the configurations of the random interlacements on the sets~$\sK_1$ and~$\sK_2$.

{\cor \label{Coro1} Suppose that we are given two functions $f_1:\{0,1\}^{\emph{\sK}_1}\to [0,1]$ and $f_2:\{0,1\}^{\emph{\sK}_2}\to [0,1]$ that depend only on the configuration of the random interlacements inside the sets $\emph{\sK}_1$ and $\emph{\sK}_2$, respectively. There exists a positive constant $C_3$, depending only on the dimension $d$, such that, if $\dist(\emph{\sK}_1,\emph{\sK}_2) \geq C_1 \diam(\emph{\sK}_1)$,
\begin{align}
\label{bound_12}
|\emph{Cov}(f_1(\I^u_{\emph{\sK}_1}), f_2(\I^u_{\emph{\sK}_2}))| 
\leq C_3 \sqrt{u} \frac{\capa(\emph{\sK}_1)^{\frac{3}{2}}}{\dist(\emph{\sK}_1,\emph{\sK}_2)^{d-2}}.
\end{align}
}

From equation $(2.15)$ of~\cite{Szn10}, we can easily obtain that, for two functions $f_1:\{0,1\}^{\sK_1}\to [0,1]$ and $f_2:\{0,1\}^{\sK_2}\to [0,1]$, that depend only on the configuration of the random interlacements inside the sets $\sK_1$ and $\sK_2$, respectively, 
\begin{align}
\label{bound_123}
|\text{Cov}(f_1(\I^u_{\sK_1}), f_2(\I^u_{\sK_2}))| 
\leq C_4 u \frac{\capa(\sK_1)^2}{\dist(\sK_1,\sK_2)^{d-2}},
\end{align}
where $C_4$ is a positive constant depending only on the dimension (see also Lemma 2.1 of \cite{Bel12}).
Thus, observe that Corollary~\ref{Coro1} improves the above bound with respect to the exponent of $\capa(\sK_1)$.
For example, when $\sK_1=\BB_r(0)$ for some $r>1$, since $\capa(\sK_1)\asymp r^{d-2}$, \eqref{bound_12} leads to 
\begin{align*}
|\text{Cov}(f_1(\I^u_{\sK_1}), f_2(\I^u_{\sK_2}))| 
\leq C_5 \sqrt{u} \Bigg(\frac{r^{\frac{3}{2}}}{\|\hat{x}\|}\Bigg)^{d-2},
\end{align*}
while \eqref{bound_123} gives
\begin{align*}
|\text{Cov}(f_1(\I^u_{\sK_1}), f_2(\I^u_{\sK_2}))| 
\leq C_6 u \Bigg(\frac{r^{2}}{\|\hat{x}\|}\Bigg)^{d-2}.
\end{align*}

It is important to stress that in this paper we are 
working without \emph{sprinkling}, that is, without 
changing the value of the parameter~$u$. When dealing
with monotone (i.e., increasing or decreasing) events,
it is often convenient to change a little
bit the value of that parameter and obtain 
some quite strong decoupling inequalities,
as shown in~\cite{AP18,PT15}. This option is unfortunately
not available for arbitrary events; as mentioned above,
\eqref{bound_123} (on which this paper improves)
was the only inequality available so far.
Among potential applications, we can mention
random walks on interlacement clusters
(such a model was considered in~\cite{FP18}) --- 
indeed, when we change the topology of the cluster 
by adding/deleting trajectories, it is not at all 
a priori clear which effect this can make on the properties
of such a random walk (and, in particular, it is unclear
why that effect should be \emph{monotone} in any
possible sense).
Also, it is naturally interesting to ask oneself
what is the \emph{correct} order of correlations
between what happens in two distant sets, for the random
interlacement model. Observe that, for a (closely related)
model of Gaussian Free Field, this order is known~\cite{PR15},
so it would be natural to conjecture that
the same also holds for the random interlacements;
in this paper, although we do not prove this conjecture,
we still do advance in its direction.

The paper is organized in the following way.
In Section~\ref{notations}, we introduce the notations used throughout the text and some basic definitions regarding the random interlacements. In Section~\ref{SimRI}, we present the construction of the excursions of the trajectories of the random interlacements between the boundaries of the sets $\sK_1\cup\sK_2$ and $\BB_R^1\cup \BB_R^2$. In the same section, we define and construct the noodle soup process. Both constructions use the technique of soft local times. This technique is also used in Section~\ref{coupling} to construct a coupling between the RI and NS processes. The probability of the complement of the corresponding coupling event will provide an upper bound for the total variation distance between the laws of the RI and NS processes. Some auxiliary results are proved in Section~\ref{toolbox}. Finally, Theorem~\ref{Maintheo} and Corollary~\ref{Coro1} are respectively proved in Sections~\ref{Proofmaintheo} and~\ref{ProofCoro1}.

\section{Notations and definitions}\label{notations}

Throughout the text, we use small~$c_1,c_2,\dots$ to denote
 global constants that appear in the results, 
and~$\gamma_1,\gamma_2,\dots$ to denote ``local constants'' that
 appear locally in the proofs, restarting the enumeration at each 
new proof.

For two functions~$f$ and~$g$, we write~$f(x)\asymp g(x)$ to 
denote that there exist positive constants~$c'$ and~$c$ depending 
only on the dimension~$d$ such that~$c'g(x)\leq f(x) \leq cg(x)$.
We will also use the convention that~$\sum_{j=k}^{i} = 0$ if~$i<k$.

In the rest of this paper, we will denote by $\sV_{R}$ and $\sK$ the sets $\BB_{R}^1\cup \BB_{R}^2$ and $\sK_1\cup\sK_2$, respectively.

For an arbitrary set~$B\subset\Z^d$, we denote its internal boundary by~$\partial B$, which is defined by
\begin{align*}
\partial B = \{x\in B : \|x-y\|=1, \text{ for some } y\in B^c\},
\end{align*} 
and we also consider the external boundary~$\partial_e B$ of~$B$, defined by
\begin{align*}
\partial_e B = \{x\in B^c : \|x-y\|=1, \text{ for some } y\in B\}.
\end{align*}

We now recall the general definition of the random interlacements process in~$\Z^d$, with~$d\geq 3$, introduced in~\cite{Szn10} (see also~\cite{CA12} and \cite{DRS14}).
This process can be viewed as an infinite random cloud of doubly-infinite simple random walk trajectories modulo time-shift with attached non-negative labels.

Formally, the random interlacements process is defined through a particular Poisson point process on a properly defined space.
To make this more precise, we begin by considering the following spaces of trajectories in~$\Z^d$,
\begin{align*}
W = &\Big\{w:\Z\rightarrow\Z^d ;\, \|w(n+1)-w(n)\| = 1 
\text{ for all }n\in\Z, \\
&~~~~~\text{ and the set } \{n:w(n)=y\} \mbox{ is finite} 
\text{ for all }y\in\Z^d \Big\}, 
\end{align*}
and
\begin{align*}
W_+ = &\Big\{w:\IN\rightarrow\Z^d ; \|w(n+1)-w(n)\| = 1
\text{ for all }n\in\IN, \\
&~~~~~\text{ and the set } \{n:w(n)=y\} \mbox{ is finite} 
\text{ for all }y\in\Z^d \Big\}, 
\end{align*}
respectively endowed with the~$\sigma$-algebras~$\W$ and~$\W_+$ generated by their respective canonical coordinates, ~$(X_n)_{n\in\Z}$ and~$(X_n)_{n\in\IN}$, and then we consider the quotient space  
\begin{align*}
W^* = W/\sim, \mbox{ where } w\sim w' \mbox{ if } w(\cdot)=w'(\cdot+k) \mbox{ for some } k\in\Z, 
\end{align*}
endowed with the~$\sigma$-algebra~$\W^*$ given by
\begin{align*}
\W^* =  \Big\{U\subset W^*: (\pi^*)^{-1}(U) \in \W \Big\},
\end{align*}
where~$\pi^*$ denotes the canonical projection from~$W$ to~$W^*$.

 Additionally, for a finite set~$B\subset\Z^d$, we also introduce the set of trajectories in~$W$ that visit~$B$:
\begin{align*}
W_B = \Big\{w\in W: X_n(w) \in B \mbox{ for some } n\in\Z \Big\}.
\end{align*}

Let~$\tau_B$ be the hitting time of a finite set~$B\subset\Z^d$,
\begin{align*}
\tau_B(w) = \inf\{n\geq 1 : X_n(w)\in B\}, \text{ for } w\in W_+,
\end{align*}
and let~$H_B$ be the entrance time of~$B$,
\begin{align*}
H_B(w) = \inf\{n\geq 0 : X_n(w)\in B\}, \text{ for } w\in W_+.
\end{align*}
For finite~$B\subset\Z^d$, we define the equilibrium measure of~$B$,
\begin{align*}
e_B(x) = P_x[\tau_B=\infty]\ind_B(x), \text{ for } x\in\Z^d,
\end{align*}
where~$P_x$ is the law of a simple random walk starting at~$x$ on~$(W_+,\W_+)$, and~$\ind_B$ is the indicator function on the set~$B$. Then the capacity of~$B$ is defined as
\begin{align*}
\capa(B) = \sum_{x\in\Z^d} e_B(x),
\end{align*}
and the harmonic measure on~$B$ is just given by~$\bar{e}_B(x)=e_B(x)/\capa(B)$, for~$x\in\Z^d$.

The \textit{random interlacements process} comes as a Poisson point process on~$(W^*\times\R_+,\W^*\otimes\B(\R_+))$ with intensity measure~$\varsigma\otimes dt$, where~$dt$ is the Lebesgue measure on~$\R_+$, $\B(\R_+)$ is the Borel $\sigma$-algebra on $\R_+$, and~$\varsigma$ is the unique~$\sigma$-finite measure on~$(W^*,\W^*)$ such that
\begin{align*}
\varsigma(\cdot\cap W_B^*) = \pi^*\circ Q_B(\cdot) , \mbox{ for any finite set } B\subset\Z^d,
\end{align*}
where~$W_B^*=\pi^*(W_B)$ and~$Q_B$ (for finite~$B\subset\Z^d$) is the finite measure on~$(W,\W)$ such that, for~$A_1,A_2\in\W_+$ and~$x\in\Z^d$,
\begin{align*}
Q_B\Big[(X_{-n})_{n\geq 0} \in A_1, X_0=x, (X_n)_{n\geq 0} \in A_2\Big] = P_x[A_1 | \tau_B=\infty] e_B(x) P_x[A_2],
\end{align*}
see~\cite{Szn10}, Theorem 1.1.

To conclude the description, consider also the space of point measures  
\begin{align*}
\Omega^* &= \Big\{ \ome=\sum_{i\geq 1}\bdelta_{(w_i^*,u_i)} : w_i^*\in W^*, u_i\in\R_+, \mbox{ and } \ome(W_B^*\times [0,u]) <\infty \\
&\hspace{5.5cm} \mbox{ for every finite set } B\subset\Z^d \mbox{ and } u\geq 0 \Big\},
\end{align*}
endowed with the~$\sigma$-algebra~$\A^*$ generated by the evaluation maps~$\ome\mapsto\ome(D)$ for~$D\in\W^*\otimes\B(\R_+)$. It is on~$(\Omega^*,\A^*)$ that the law of the above mentioned Poisson point process is usually considered. 

Finally, for~$\ome\in\Omega^*$, the \textit{random interlacements} at level~$u$ is defined to be the (random) set
\begin{align*}
\I^u(\ome) = \bigcup_{i:u_i\leq u} \mbox{Range}(w_i^*),
\end{align*}
so that the trace left on a finite set~$B\subset\Z^d$ by the random interlacements at level~$u$ is just~$\I^u(\ome)\cap B$. For the sake of brevity, in the rest of this paper we will denote this last random set by~$\I^u_B$.

\section{Constructions using soft local times}
\label{SimRI} 

The main purpose of this section is to show how the soft local times method can be used to construct the random interlacements (restricted to $\sK$) and the NS processes. Besides, these constructions will be useful in Section \ref{coupling}, where we construct a coupling between both processes. We refer to \cite{PT15} for a presentation of the soft local times method.

\begin{figure}[h]
	\centering \includegraphics{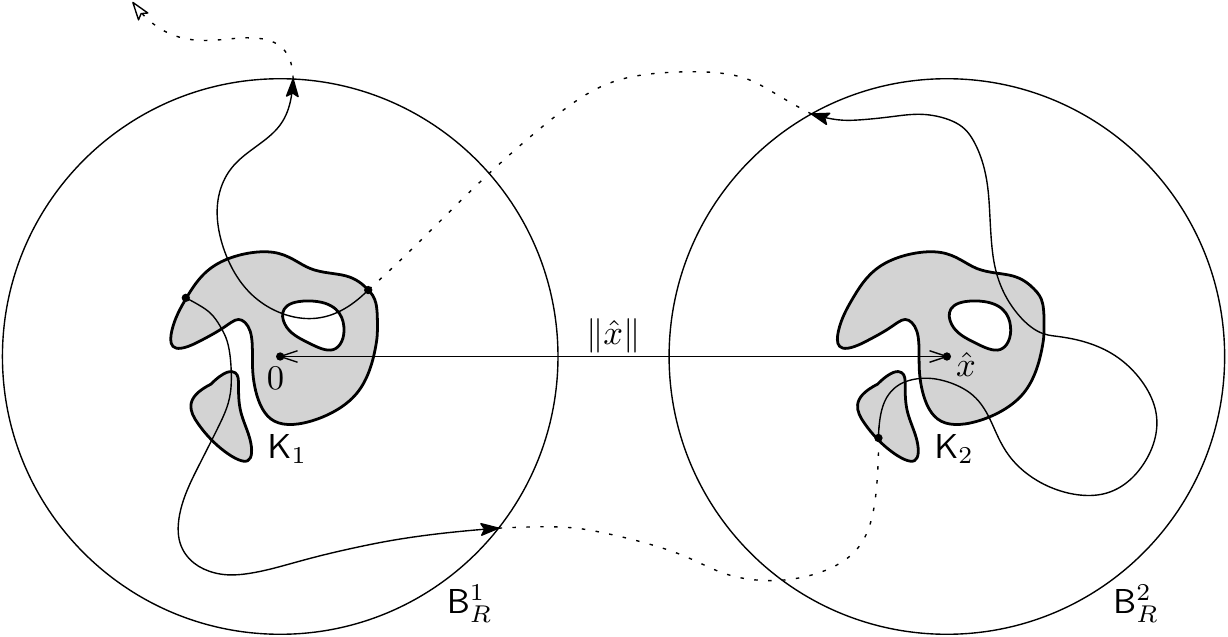} 
	\caption{The excursions of a SRW trajectory between $\partial\sK$ and $\partial_e\sV_R$.}
	\label{balls_fig}
\end{figure}

Recall that we denote by $\sK$ the set $\sK_1\cup\sK_2$, and by $\sV_R$ the set $\BB_{R}^1\cup \BB_{R}^2$. 
We suppose that $\|\hat{x}\|~\geq 4\diam(\sK_1)+3$ and we take 
$R=\frac{\|\hat{x}\|-1}{2}$. Observe that, with this choice of~$R$,
the sets $\partial_e \BB_R^1$ and $\partial_e \BB_R^2$ are 
disjoint and $\sK_1\subset \BB_R^1$ (see Figure~\ref{balls_fig} for illustration).

\subsection{Construction of random interlacements restricted to~$\sK$} \label{SimRI_2}

The random interlacements (at level~$u$) restricted to some finite subset of~$\Z^d$, say~$B$, can be obtained as the trace on~$B$ of a Poisson number, with parameter~$u\capa(B)$, of independent simple 
random walk trajectories starting at the internal boundary~$\partial B$. The starting sites of each such trajectory are independently chosen according to the harmonic measure on~$B$.

This description will allow us to construct the random interlacements restricted to~$\sK$ from a soup of excursions between~$\partial \sK$ and~$\partial_e \sV_R$.

We denote by~$\II^u$ the random multiple set of excursions 
(of the trajectories of the random interlacements process at 
level~$u$) between~$\partial \sK$ and~$\partial_e \sV_R$. 
Thus, $\I^u_{\sK}$ will be the trace left on~$\sK$ by the
 elements of~$\II^u$.
Next we present a method to construct~$\II^u$, in the same spirit 
as~\cite{BGP16}.

First, we introduce the set of simple random walk excursions between~$\partial \sK$ and~$\partial_e \sV_R$,
\begin{align*}
\Sigma = \Big\{ z=(x_0,x_1,\dots,x_n) : n<\infty,&~ \|x_{i-1}-x_{i}\|=1, i=1,2,\dots,n, \\
& x_0\in \partial \sK, \mbox{ and } x_j\in\partial_e \sV_R \text{ only for } j=n \Big\}.
\end{align*} 
That is, $\Sigma$ is the set of finite nearest neighbour paths on~$\Z^d$, starting at~$\partial \sK_1$ or~$\partial \sK_2$ and ending at their first visit to~$\partial_e \BB^1_R$ or~$\partial_e \BB^2_R$. 
Note that~$\partial_e \sV_R$ separates the sets $\sK_1$ and $\sK_2$, in the sense that any trajectory that goes from one of them to the other must cross~$\partial_e \sV_R$.

Observe that, in the interlacements process restricted to~$\sK$, 
each one of the simple random walk trajectories starting 
at~$\partial \sK$ will visit~$\partial_e \sV_R$ at a finite time, 
thus performing an excursion (which belongs to~$\Sigma$), and from this moment it 
can return to one of the sets $\sK_1$ or $\sK_2$ and then perform another 
excursion, or it can never return to them. Moreover, any 
trajectory makes only a finite number of such excursions (see Figure \ref{balls_fig}).

Now, we consider the infimum of the probabilities of escaping the set $\sK$, starting at~$\partial_e \sV_R$, namely
\begin{align*}
q = q(R, \capa(\sK)) := \inf_{y\in \partial_e \sV_R} P_y[\tau_{\sK}=\infty].
\end{align*}
Then we consider $N_1$ and $N_2$, independent Poisson random variables with parameters respectively equal to 
\begin{align*}
(1-q) u\capa(\sK) ~~\text{  and  }~~  q u\capa(\sK).
\end{align*}

One way to construct the excursions performed by the simple 
random walks of the random interlacements process 
between~$\partial \sK$ and~$\partial_e \sV_R$ is through 
the technique of soft local times, recently introduced 
in~\cite{PT15} (see also~\cite{CGPV13}).

For our purposes in this paper, we will use the soft local times to construct~$\II^u$ in a slightly different way compared to the original approach of~\cite{PT15}, by obtaining first~$N_1$ ``possibly returning'' trajectories and then~$N_2$ ``non-returning'' trajectories of the random interlacements process.

We consider a Poisson point process~$\eta=\sum_{\lambda\in\Lambda} \bdelta_{(z_{\lambda},t_{\lambda})}$ on~$\Sigma\times\R_+$ (here~$\Lambda$ is a countable set) with intensity measure given by~$\mu\otimes dt$, where~$dt$ is the Lebesgue measure on~$\R_+$ and~$\mu$ is the measure on~$\Sigma$ defined by
\begin{align*}
\mu(\Y) = \sum_{x\in\partial \sK} \bar{e}_{\sK}(x) P_x[(X_0, X_1,\dots,X_{H_{\partial_e \sV_R}})\in\Y], \text{ for } \Y\in\mS,
\end{align*}
with~$\mS$ being the~$\sigma$-algebra generated by the canonical coordinates.

Before proceeding, we need to introduce more notations. 
Let us define the map~$S_0$ from the space of excursions
to the boundary of the set $\sK$, 
$S_0: \Sigma \rightarrow \partial \sK$, 
which selects the starting point of each excursion, 
\begin{align*}
S_0(z)= x_0, \text{ for } z=(x_0,x_1,\dots,x_n) \in \Sigma,
\end{align*}
and also the map~$S_f$ from the space of excursions to the external boundary of the balls $\BB_R^1$ and $\BB_R^2$, $S_f: \Sigma \rightarrow \partial_e \sV_R$, which selects the last point of each excursion, 
\begin{align*}
S_f(z)= x_n, \text{ for } z=(x_0,x_1,\dots,x_n) \in \Sigma.
\end{align*}
Further, we define
\begin{align*}
p_{S_f(z)} := P_{S_f(z)}[\tau_{\sK}=\infty],
\end{align*}
for any excursion~$z\in\Sigma$.

We start with the construction of the excursions of the first~$N_1$ trajectories of the random interlacements process. For that, we consider a family~$(\zeta_i)_{i\geq 1}$ of independent~$[0,1]$-uniformly distributed random variables, which is also independent of all the other random elements. Then, using the point process~$\eta$ described above and~$N_1$, we proceed with the soft local times scheme as follows.

For the construction of the successive excursions of the first random walk trajectory, we start constructing its first excursion, which we call~$z_1$, by defining
\begin{align*}
\xi_1 &= \inf\Big\{\ell\geq 0: \text{ there exists } \lambda
 \text{ such that } \ell \geq t_{\lambda}\Big\}, \\
G^{\II}_1(z) &=  \xi_1,
\end{align*}
and~$(z_1,t_1)$ to be the unique pair~$(z_{\lambda},t_{\lambda})$ satisfying~$G^{\II}_1(z_{\lambda}) = t_{\lambda}$.

Then, once we have obtained the first excursion~$z_1$ of the first trajectory, we decide whether the same trajectory will perform another excursion or not. For that, we consider the family~$(\zeta_i)_{i\geq 1}$ and we use the following criterion:
\begin{itemize}
\item If~$\zeta_1 \leq (p_{S_f(z_1)}-q)/(1-q)$ (which is equivalent to a uniform r.v.~in $[q,1]$ be smaller than $p_{S_f(z_1)}$), then we stop the construction of the excursions of this trajectory, and we proceed with the construction of the second trajectory (we interpret this situation by saying that this is a non-returning trajectory, which therefore makes only its first excursion and then escapes to infinity).
\item Otherwise, we keep constructing the subsequent excursions of the same trajectory, as we describe below. 
\end{itemize}

Define the transition density~$g$ (with respect to the measure~$\mu$)
\begin{align}
\label{trans_g}
g(z,z') = \frac{P_{S_f(z)}[X_{\tau_{\sK}}=S_0(z')|\tau_{\sK}<\infty]}{\bar{e}_{\sK} ( S_0(z'))},
\end{align}
and for~$n=2,3,\dots$, define
\begin{align*}
\xi_n &= \inf\Big\{\ell\geq 0:  \text{ there exists }  (z_{\lambda},t_{\lambda})\notin\{(z_k,t_k)\}_{k=1}^{n-1} \text{ such that } \\  & \hspace{5cm}G^{\II}_{n-1}(z_{\lambda}) + \ell g(z_{n-1},z_{\lambda}) \geq t_{\lambda}\Big\}, \\
G^{\II}_n(z) &= G^{\II}_{n-1}(z) + \xi_n  g(z_{n-1},z),
\end{align*}
and~$(z_n,t_n)$ to be the unique pair~$(z_{\lambda},t_{\lambda})$ out of the set~$\{(z_k,t_k)\}_{k=1}^{n-1}$ satisfying~$G^{\II}_n(z_{\lambda}) = t_{\lambda}$. 
For each value of~$n$, at the end of the corresponding step we obtain the $n$-th excursion~$z_n$ of the first trajectory, and we decide whether the trajectory will perform another excursion or not analogously to the first step, but now with a slightly different criterion involving the corresponding uniform random variable: 
\begin{itemize}
\item If~$\zeta_n \leq p_{S_f(z_n)}$, then we stop the construction of the excursions of this trajectory, and we proceed to the construction of the second trajectory.
\item Otherwise, we keep constructing the subsequent excursions of the same trajectory, iteratively in~$n$. 
\end{itemize}

Thus, we will conclude the construction of the first trajectory at the random time 
\begin{align*}	
T_1 = \inf\{n\geq 1: \zeta_n \leq \kappa_n\},
\end{align*}
where 
\begin{align*}	
\kappa_1 &= \frac{p_{S_f(z_1)}-q}{1-q}, \\ 
\text{ and } \kappa_n &= p_{S_f(z_n)}, \text{ for } n=2,3,\dots, T_1,
\end{align*}
and at this moment we obtain the accumulated soft local time corresponding to this first trajectory,
\begin{align*}
G^{\II}_{T_1}(z) = \sum_{j=1}^{T_1} \xi_j g(z_{j-1},z),
\end{align*}
for~$z\in\Sigma$, with the convention that~$g(z_0,z) = 1$.

Then we proceed with the construction of the remaining~$N_1-1$ trajectories, just as we would do in the original soft local times procedure, using the same point process~$\eta$, but now imitating the above described iterative scheme to construct the excursions of each remaining trajectory. That is, we begin with the density~$\ind_{\Sigma}$ for the construction of the first excursion of a new trajectory and then we continue with the transition density~$g$ for the construction of the next excursions, implementing the comparisons involving the uniform random variables~$(\zeta_i)_{i\geq 1}$ to decide when to finish the construction of each trajectory. More precisely, when dealing with the the~$j$-th trajectory ($ j\geq 2$), we use the random variables~$\zeta_k$ with  
\begin{align}
\label{k_values}
k= \Theta_{j-1} +1, \Theta_{j-1} +2, \dots, \Theta_{j}, ~\text{ where }~ \Theta_{j} = \sum_{i=1}^{j \wedge N_1}T_i,
\end{align}
(using the convention that $\Theta_1=T_1$) to decide when to finish the construction of this trajectory. This procedure generates the corresponding families~$(\xi_k)_{k}$,  $(G^{\II}_k(\cdot))_{k}$ and~$(z_k)_{k}$ for the values of~$k$ as in~\eqref{k_values}, and lasts the corresponding random time 
\begin{align*}	
T_j = \inf\Big\{n\geq 1: \zeta_{\Theta_{j-1} +n} \leq \kappa_{\Theta_{j-1} +n}\Big\},
\end{align*}
where 
\begin{align*}	
\kappa_{\Theta_{j-1} +1} &= \frac{p_{S_f(z_{\Theta_{j-1} +1})}-q}{1-q}, \\
\text{ and } \kappa_n &= p_{S_f(z_n)}, \text{ for } n=\Theta_{j-1} +2, \Theta_{j-1} +3, \dots, \Theta_{j}.
\end{align*}

At the end of this procedure, we obtain the excursions of the first~$N_1$ trajectories of the random interlacements at level~$u$, $(z_k : 1\leq k\leq \Theta_{N_1})$, and we also obtain the accumulated soft local time up to the~$N_1$-th trajectory, that is, the accumulated soft local time corresponding to the possibly returning excursions, 
\begin{align*}
 G^{\II}_{\Theta_{N_1}}(z) = \sum_{k=1}^{\Theta_{N_1}} \xi_k \tilde{g}(z_{k-1},z), ~\text{ for } z\in\Sigma,
\end{align*} 
with the convention that $\Theta_0=0$ and
\begin{align*}
\tilde{g}(z_{k-1},z) = 
\begin{cases} 
1  & \mbox{ for } k-1=\Theta_{j}, j=0,1,2, \dots, N_1-1, \\   
g(z_{k-1},z)& \mbox{ otherwise. }
\end{cases} 
\end{align*}

Observe that, since the transition density~$g(\cdot,z)$ depends on~$z$ only through its initial point, the same happens to the accumulated soft local time~$G^{\II}_{\Theta_{N_1}}(z)$. Also, note that the random variables~$(T_j)_{j\geq 1}$ are all independent and identically distributed.

Finally, we complete~$\II^u$ by obtaining the remaining non-returning trajectories. We use~$N_2$ and the points of the Poisson point process~$\eta$ left above the curve~$G^{\II}_{\Theta_{N_1}}$ after the previously described construction. Then, for~$n=\Theta_{N_1} +j$, $j=1,2,\dots,N_2$, define
\begin{align*}
\xi_n &= \inf\Big\{\ell\geq 0: \exists (z_{\lambda},t_{\lambda})\notin\{(z_k,t_k)\}_{k=1}^{n-1} \text{ such that }  G^{\II}_{n-1}(z_{\lambda}) + \ell \geq t_{\lambda}\Big\}, \\
G^{\II}_n(z) &= G^{\II}_{n-1}(z) + \xi_n ,
\end{align*}
and~$(z_n,t_n)$ to be the unique pair~$(z_{\lambda},t_{\lambda})$ out of the set~$\{(z_k,t_k)\}_{k=1}^{n-1}$ satisfying~$G^{\II}_n(z_{\lambda}) = t_{\lambda}$.
Thus, introducing~$\N:=\Theta_{N_1}+N_2$, at the end of these iterations we obtain the accumulated soft local time~$G^{\II}_{\N}$ corresponding to the process~$\II^u$,
\begin{align*}
G^{\II}_{\N}(z) = G^{\II}_{\Theta_{N_1}}(z) +  \sum_{k=\Theta_{N_1}+1}^{\N} \xi_k,
\end{align*}
for~$z\in\Sigma$.

In Section~\ref{coupling} we present a slightly different construction using the soft local times, in order to couple~$\II^u$ with the noodle soup process which we describe in the next section.

\subsection{Definition and construction of the NS process}\label{MES_process}

Let us now describe the noodle soup process, 
at level~$u$, which will be denoted by~$\MM^u$. 

First, recall the random variables~$N_1$, $N_2$, $\{T_j\}_{j=1}^{N_1}$, $\Theta_{N_1}$ and~$\N$ from the last section.
From the construction of that section, for each~$j\in\{1,2,\dots,N_1\}$, observe that the random quantity~$T_j$ represents the (random) number of excursions that is performed by the~$j$-th trajectory, so that the total number of excursions performed by all the possibly returning trajectories in that construction is just~$\Theta_{N_1}$.

Next, consider a family~$(\E_j)_{j\geq 1}$, of independent simple random walk excursions in~$\Sigma$, each excursion starting according to the harmonic measure~$\bar{e}_{\sK}(\cdot)$, and then running up to its first visit to the separating set~$\partial_e \sV_R$. Also, consider a Poisson random variable~$\N'$, independent of~$(\E_j)_{j\geq 1}$, with the same mean as~$\N$. The NS process~$\MM^u$ is simply defined as the multiple set of excursions~$\E_1, \E_2, \dots, \E_{\N'}$, and we denote by~$\M^u_{\sK}$ the trace left on~$\sK$ by the excursions of~$\MM^u$.

Thus, by definition, the process~$\MM^u$ has the same expected number of random walk excursions as in~$\II^u$, but in~$\MM^u$ all the excursions are completely independent, unlike what happens in $\II^u$, which has a ``portion'' of dependent excursions. 

Before ending this section, we just remark that it is possible to apply the soft local times technique to construct~$\MM^u$, as we briefly describe now. 
We use~$\N'$ and the Poisson point process~$\eta=\sum_{\lambda\in\Lambda} \bdelta_{(z_{\lambda},t_{\lambda})}$ on~$\Sigma\times\R_+$ with intensity measure given by~$\mu\otimes dt$, introduced in the last section. Define~$G^{\MM}_0(\cdot)\equiv 0$, and for~$n=1,2,\dots,\N'$,
\begin{align*}
\xi'_n &= \inf\Big\{\ell\geq 0: \exists (z_{\lambda},t_{\lambda})\notin\{(z_k,t_k)\}_{k=1}^{n-1} \text{ such that } G^{\MM}_{n-1}(z_{\lambda}) + \ell \geq t_{\lambda}\Big\}, \\
G^{\MM}_n(z) &= G^{\MM}_{n-1}(z) + \xi'_n,
\end{align*}
and~$(z_n,t_n)$ to be the unique pair~$(z_{\lambda},t_{\lambda})$ out of the set~$\{(z_k,t_k)\}_{k=1}^{n-1}$ satisfying~$G^{\MM}_n(z_{\lambda}) = t_{\lambda}$.
In this way we obtain the accumulated soft local time~$G^{\MM}_{\N'}$ corresponding to the process~$\MM^u$, namely
\begin{align*}
G^{\MM}_{\N'}(z) =  \sum_{k=1}^{\N'} \xi'_k,
\end{align*}
for~$z\in\Sigma$.

\section{Coupling between RI and NS processes}\label{coupling}

In this section we present the construction of a coupling between $\II^u$ and $\MM^u$, using the soft local times technique. This coupling is inspired from the one described in \cite{BGP16}, Section 4.

We will construct two copies of the Poisson point process~$\eta$ on $\Sigma\times\R_+$ with intensity $\mu\otimes dt$, which we call~$\eta^{\II^u}$ and~$\eta^{\MM^u}$. These copies will be such that the construction of Section~\ref{SimRI_2} applied to~$\eta^{\II^u}$ and the construction of Section~\ref{MES_process} applied to~$\eta^{\MM^u}$  will give high probability of successful coupling between $\MM^u$ and $\II^u$, when $\|\hat{x}\|$ is large.

\begin{figure}[]
	\centering \includegraphics[scale=0.8]{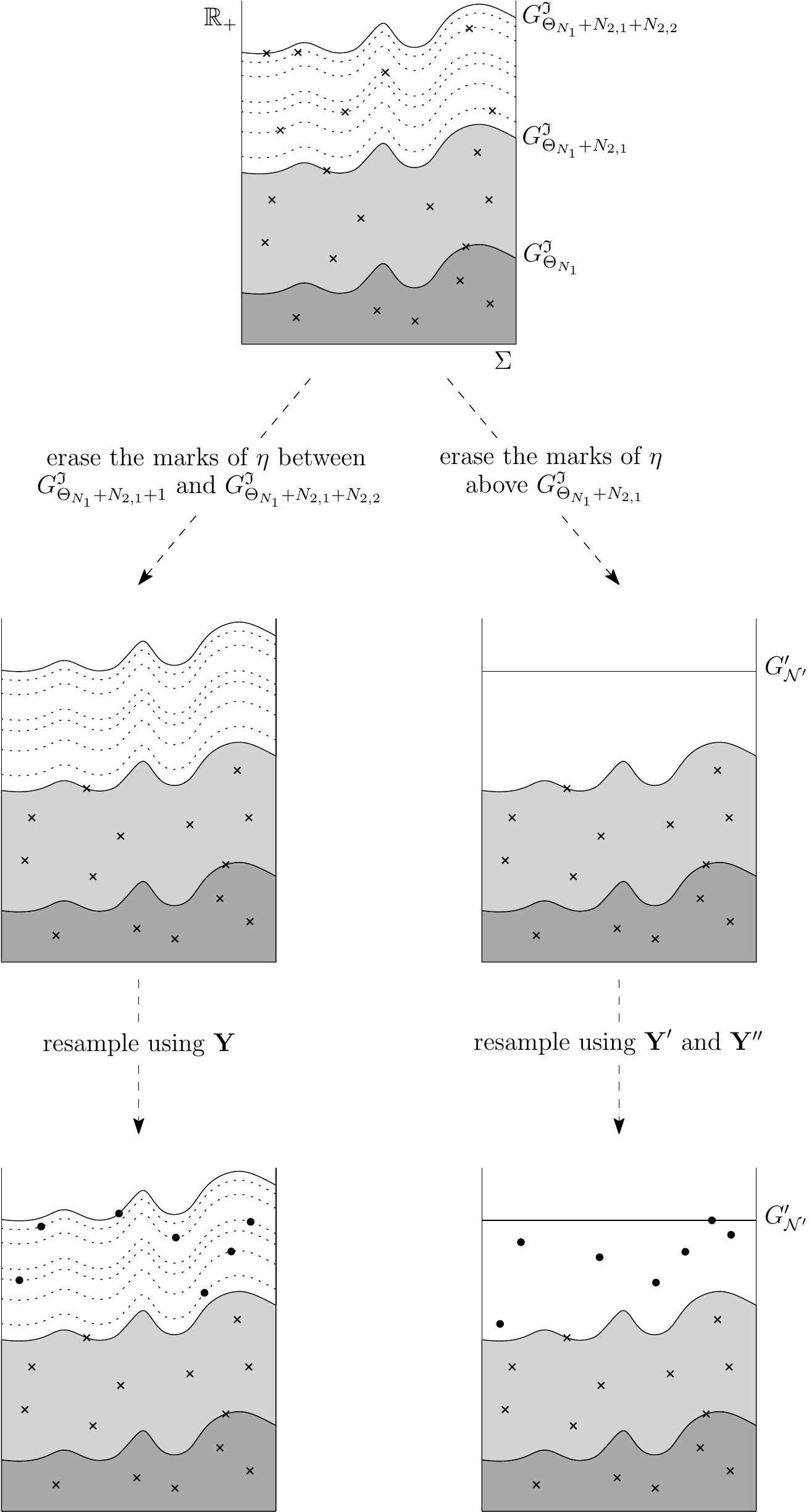} 
	\caption{Construction of the coupling.}
	\label{coupling_fig}
\end{figure}

At this point observe that, differently of what happens in \cite{BGP16}, where we have two processes of the same size $n\in\IN$, here we have two processes ($\II^u$ and $\MM^u$) with random cardinalities ($\N$ and $\N'$, respectively) which have different laws. Hence, in a first step we will couple this random cardinalities (see \eqref{eventD}), and then we will use a ``resampling'' technique (like in \cite{BGP16}) to complete the construction of the coupling between $\II^u$ and $\MM^u$ (see Figure~\ref{coupling_fig}).

Fix once for all an enumeration of the sites of $\partial \sK$. 
This will allow us to consider vectors of the form 
$(z_j)_{j\in \partial \sK}$ without any ambiguity. Also, suppose that on some probability space $(\Omega, \mathcal{F}, \IP)$, we are given the following independent random elements:

\begin{itemize}
\label{families}
\item $N_1$, $N'_1$ and~$N_{2,2}$, independent Poisson r.v.\ with respective parameters~$(1-q)u\capa(\sK)$, $\IE[\Theta_{N_1}]$ and~$\frac{q}{2}u\capa(\sK)$;
\item a sequence~$(\zeta_i)_{i\geq 1}$ of independent Uniform(0,1) r.v.;
\item a Poisson point process $\eta$ on $\Sigma\times\R_+$ with intensity $\mu\otimes dt$;
\end{itemize}
The random elements  $N_1$, $(\zeta_i)_{i\geq 1}$ and $\eta$ will have the same role as in Section \ref{SimRI_2} for the construction of $\II^u$ and  $N'_1$, $N_{2,2}$ will be used as part of the number of excursions of $\MM^u$.
Further up, some other random elements will be defined, and we assume that $(\Omega, \mathcal{F}, \IP)$ is large enough to support all the random elements to be defined in this section.

We start with the construction of~$\eta^{\II^u}$. For that, we use the random variables~$(\zeta_i)_{j\geq 1}$, $N_1$ and the Poisson point process $\eta$ to construct the first portion of excursions of~$\II^u$ through the soft local times procedure exactly as described in Section~\ref{SimRI_2}. This step will produce~$\Theta_{N_1}$ excursions, corresponding to the projections of the points of~$\eta$ below~$G^{\II}_{\Theta_{N_1}}$, on~$\Sigma$. In this way, we obtain the soft local times curves~$G^{\II}_{i}$, $1\leq i\leq \Theta_{N_1}$, together with the sequences~$\xi_1, \xi_2,\dots, \xi_{\Theta_{N_1}}$ and~$z_1, z_2,\dots, z_{\Theta_{N_1}}$.

Now, we continue with the construction of the coupling between the number of excursions of $\II^u$ and $\MM^u$. One difficulty arises here since the random variable $\Theta_{N_1}$ depends on some excursions of $\II^u$. To cope with this problem, we proceed as follows. Consider two families of random variables $(\tilde{X}_k)_{k\in\Z}$ and~$(\tilde{Y}_k)_{k\in\Z}$ such that~$((\tilde{X}_k)_{k\in\Z},(\tilde{Y}_k)_{k\in\Z})$ is independent of all the other random elements already introduced. The elements composing each family are all independent Poisson random variables with parameter~$\frac{q}{2}u\capa(\sK)$ and, for each~$k\in\Z$, $\tilde{Y}_k$ and~$k+\tilde{X}_k$ are maximally coupled. Then we define the random variables~$N_{2,1}$ and~$N'_{2,1}$ such that, for all~$k\in\Z$,  
\begin{align}
\label{def_N21}
N_{2,1} = \tilde{X}_k  ~\mbox{ and }~ N'_{2,1} = \tilde{Y}_k
\end{align}
on the event~$\{\Theta_{N_1}-N'_1 = k\}$. Note that $(\Theta_{N_1},N_1)$, $N_{2,1}$ and $N_{2,2}$ are independent and $N_{2,1}+N_{2,2}$ has the same law as~$N_2$ from Section~\ref{SimRI_2}. We will use $\Theta_{N_1}+N_{2,1}+N_{2,2}$ as the number of excursions for $\II^u$. Also, observe that $N'_1$, $N'_{2,1}$ and $N_{2,2}$ are independent and $N'_1+N'_{2,1}+N_{2,2}$ has the same law as $\N'$ from Section \ref{MES_process}. We will use $N'_1+N'_{2,1}+N_{2,2}$ as the number of excursions for $\MM^u$.

Next, we use the Poisson point process $\eta$ to complete the construction of~$\II^u$, again as described in Section~\ref{SimRI_2}, but now using~$N_{2,1} + N_{2,2}$ instead of~$N_2$. In this way, we obtain the soft local times curves~$G^{\II}_{i}$, $\Theta_{N_1}+1\leq i\leq \Theta_{N_1}+N_{2,1} + N_{2,2}$ together with the sequences
\begin{align*}
\xi_{\Theta_{N_1}+1}, \xi_{\Theta_{N_1}+2},\dots, \xi_{\Theta_{N_1}+N_{2,1} + N_{2,2}} ~\text{ and }~ z_{\Theta_{N_1}+1}, z_{\Theta_{N_1}+2},\dots, z_{\Theta_{N_1}+N_{2,1} + N_{2,2}}.
\end{align*}

Then, we use a ``resampling'' scheme: we first ``erase'' all marks of 
the point process~$\eta$ that are on the curves
 $G^{\II}_{\Theta_{N_1}+N_{2,1} +1},
\dots,G^{\II}_{\Theta_{N_1}+N_{2,1} +N_{2,2}}$, 
and then reconstruct the process in the following way.
Consider the random element 
\[
W:=\Big(N_1, N'_1, N_{2,1}, N'_{2,1}, N_{2,2}, (\zeta_i)_{i\geq 1}, \eta\Big)
\]
and introduce the random vector ${\bf Y}:=(Y_1,\dots,Y_{N_{2,2}})$ such that, under $\IP[\;\cdot \mid W=w]$, for $w\in\{N_{2,2}\geq 1\}$, its coordinates are independent simple random walk excursions (of the space~$\Sigma$) with initial law given by the harmonic measure $\bar{e}_{\sK}$.
We use the random variables $\xi_{\Theta_{N_1}+N_{2,1} +1}, \xi_{\Theta_{N_1}+N_{2,1} +2}, \dots,\xi_{\Theta_{N_1}+N_{2,1} +N_{2,2}}$ and the vector ${\bf Y}$ to reconstruct~$\eta$ from $G^{\II}_{\Theta_{N_1}+N_{2,1} +1}$ to $G^{\II}_{\Theta_{N_1}+N_{2,1} +N_{2,2}}$, by placing the marks
\begin{align*}
\Big(Y_1,G^{\II}_{\Theta_{N_1}+N_{2,1} +1}(Y_1)\Big),
&\Big(Y_2,G^{\II}_{\Theta_{N_1}+N_{2,1} +2}(Y_2)\Big),\dots \\ 
&\dots, \Big(Y_{N_{2,2}},G^{\II}_{\Theta_{N_1}+N_{2,1} +N_{2,2}}(Y_{N_{2,2}})\Big),
\end{align*}
on the curves $G^{\II}_{\Theta_{N_1}+N_{2,1} +1},\dots,G^{\II}_{\Theta_{N_1}+N_{2,1} +N_{2,2}}$ (see Figure \ref{coupling_fig}). Observe that this step is immaterial if~$N_{2,2}$ is equal to zero. Hence, $\eta^{\II^u}$ is the point process obtained using this resampling procedure.

Now, let us continue with the construction of~$\eta^{\MM^u}$. 
For that, we introduce the random variable
\begin{align*}
\Xi_{2,2}:=\sum_{j=\Theta_{N_1}+N_{2,1}+1}^{\Theta_{N_1}+N_{2,1}+N_{2,2}}\xi_j, 
\end{align*}
and the random function 
\begin{align*}
G'_{\N'}(z)=\inf_{z'\in\Sigma}G^{\II}_{\Theta_{N_1}+N_{2,1}}(z') + \Xi_{2,2}, ~\text{ for } z\in\Sigma.
\end{align*}

Given $G^{\II}_{\Theta_{N_1}+N_{2,1}}$ and~$G'_{\N'}$, we also consider the conditional probability law $\Psi$ such that, for~$z\in\Sigma$,
\begin{align*}
\Psi(z) = 
\displaystyle\frac{\Big(G'_{\N'}(z)-G^{\II}_{\Theta_{N_1}+N_{2,1}}(z)\Big)_+}{\displaystyle\int_{\Sigma}\Big(G'_{\N'}(z)-G^{\II}_{\Theta_{N_1}+N_{2,1}}(z)\Big)_+ \mu(dz)}
\end{align*}
if the denominator of the last expression is positive, and $\Psi(z)=\ind_{\Sigma}(z)$ otherwise. 
Anticipating on what is coming, the law $\Psi\text{d}\mu$ will be used to sample the first coordinates of the~$N_{2,2}-1$ marks of~$\eta^{\MM^u}$ in the open set of $\Sigma\times\R_+$ delimited by~$G^{\II}_{\Theta_{N_1}+N_{2,1}}$ and~$G'_{\N'}$.

We first construct the marks of~$\eta^{\MM^u}$ below~$G'_{\N'}$.
For this, consider the random vector ${\bf Y'}:=(Y'_1,\dots,Y'_{N_{2,2}})$ and the random element ${\bf Y''}:=((Y''_1(z),\dots,Y''_{N_{2,2}}(z))_{z\in\Sigma}$ such that, under $\IP[\;\cdot \mid W=w]$, for $w\in\{N_{2,2}\geq 1\}$, ${\bf Y'}$ and ${\bf Y''}$ are independent and
\begin{itemize}
\item $Y'_{N_{2,2}}$ has law $\frac{\ind{\{G'_{\N'}>G^{\II}_{\Theta_{N_1}+N_{2,1}}\}}}{\mu[G'_{\N'}>G^{\II}_{\Theta_{N_1}+N_{2,1}}]}d\mu$ if $\mu[G'_{\N'}>G^{\II}_{\Theta_{N_1}+N_{2,1}}]>0$, and law $\mu$ otherwise. Also, $Y'_{N_{2,2}}$ is maximally coupled with $Y_{N_{2,2}}$;
\item $(Y'_1,\dots,Y'_{N_{2,2}-1})$ has its coordinates being independent with law~$\Psi\text{d}\mu$, and the random elements $(\sum_{i=1}^{N_{2,2}-1}{\bdelta}_{x}(Y'_i))_{x\in \Sigma }$ and $(\sum_{i=1}^{N_{2,2}-1}{\bdelta}_{x}(Y_i))_{x\in \Sigma }$ are maximally coupled;
\item the elements $(Y''_{1}(z))_{z\in\Sigma}, \dots, (Y''_{N_{2,2}}(z))_{z\in\Sigma}$ are i.i.d.;
\item $(Y''_{1}(z))_{z\in\Sigma}$ is a family of independent random variables and, for $z\in\Sigma$, $Y''_{1}(z)$ has law $U(0, G'_{\N'}(z)-G^{\II}_{\Theta_{N_1}+N_{2,1}}(z))$ if $\Psi(z)>0$, and law $U(0,1)$ if $\Psi(z)=0$.
\end{itemize}

We construct the marks of the point 
process~$\eta^{\MM^u}$ below $G'_{\N'}$, using ${\bf Y'}$ and ${\bf Y''}$ in the following way:
we keep the marks obtained below $G^{\II}_{\Theta_{N_1}+N_{2,1}}$ 
and we use the law $\Psi\text{d}\mu$ to complete the process 
until~$G'_{\N'}$.
For this, we adopt a resampling scheme as before. 
We first erase all the marks of the point process~$\eta$
that are (strictly) above $G^{\II}_{\Theta_{N_1}+N_{2,1}}$, 
then we resample the part of the process $\eta$ up to~$G'_{\N'}$, 
using the marks:
\begin{align*}
\Big(Y'_1,& G^{\II}_{\Theta_{N_1}+N_{2,1}}(Y'_1)+ Y''_1(Y'_1)\Big),\dots, \\ 
& \Big(Y'_j,G^{\II}_{\Theta_{N_1}+N_{2,1}}(Y'_j)+Y''_j(Y'_j)\Big),\dots, \\
&~~~~~~\Big(Y'_{N_{2,2}-1},G^{\II}_{\Theta_{N_1}+N_{2,1}}(Y'_{N_{2,2}-1})+Y''_{N_{2,2}-1}(Y'_{N_{2,2}-1})\Big), \Big( Y'_{N_{2,2}}, G'_{\N'}(Y'_{N_{2,2}}) \Big)
\end{align*}
(see Figure \ref{coupling_fig}). 


Then, to complete the marks of~$\eta^{\MM^u}$ above~$G'_{\N'}$, we ``glue'' a copy of~$\eta$, independent of everything, above~$G'_{\N'}\vee G^{\II}_{\Theta_{N_1}+N_{2,1}}$.

Finally, we construct the process $\MM^u$ using the soft local times technique exactly as described in Section \ref{MES_process}, but now using the quantity $N'_1+N'_{2,1}+N_{2,2}$ instead of $\N'$, and the Poisson point process $\eta^{\MM^u}$ instead of $\eta$.

Let us denote by $\Pi(G^{\II}_{\Theta_{N_1}+N_{2,1} +N_{2,2}})$ and $\Pi(G^{\MM}_{N'_1+N'_{2,1}+N_{2,2}})$ the multiple sets formed by the projections on $\Sigma$ of the marks of~$\eta^{\II^u}$ below the curve $G^{\II}_{\Theta_{N_1}+N_{2,1} +N_{2,2}}$ and of~$\eta^{\MM^u}$ below $G^{\MM}_{N'_1+N'_{2,1}+N_{2,2}}$, respectively. 


We have the following
{\prop 
	We have that~$\eta^{\II^u}\eqlaw\eta^{\MM^u}\eqlaw\eta$. 
	Furthermore, it holds that
	$\Pi(G^{\II}_{\Theta_{N_1}+N_{2,1} +N_{2,2}})\eqlaw\II^u$ and 
	$\Pi(G^{\MM}_{N'_1+N'_{2,1}+N_{2,2}})\eqlaw\MM^u$ 
	(where $\eqlaw $ stands for equality in law).\\
}

\begin{proof}
	The proof follows from elementary properties of Poisson point processes and the fact that, by construction, $\Theta_{N_1}$ and $N_{2,1}$ are independent, and also that~$\eta^{\MM^u}$ is independent of the sum $N'_1+N'_{2,1}+N_{2,2}$. Indeed, for all nonnegative integers $k$ and $\ell$, 
	\begin{align*}
	\IP[\Theta_{N_1} = k, N_{2,1} = \ell] &= \sum_{m=0}^{\infty} \IP[\Theta_{N_1} = k, N_{2,1} = \ell, N'_1=m]\\
	&= \sum_{m=0}^{\infty} \IP[\Theta_{N_1}-N'_1 = k-m, N_{2,1} = \ell, N'_1=m]\\
	&= \sum_{m=0}^{\infty} \IP[\Theta_{N_1}-N'_1 = k-m, \tilde{X}_{k-m} = \ell, N'_1=m]\\
	&= \sum_{m=0}^{\infty} \IP[\Theta_{N_1}-N'_1 = k-m, N'_1=m] \IP[\tilde{X}_{k-m} = \ell]\\
	&= \IP[N_{2,1} = \ell]\sum_{m=0}^{\infty} \IP[\Theta_{N_1}-N'_1 = k-m, N'_1=m] \\
	&= \IP[N_{2,1} = \ell]\IP[\Theta_{N_1} = k]. \\
	\end{align*}
	Then, observe that to check that~$\eta^{\MM^u}$ is independent of the sum $N'_1+N'_{2,1}+N_{2,2}$, it is enough to check that $\Theta_{N_1}$, $N'_1$ and $N'_{2,1}$ are independent. For all nonnegative integers $k$, $\ell$ and $m$, we have
	\begin{align*}
	\IP[\Theta_{N_1} = k, N'_{2,1} = \ell, N'_1=m] &= \IP[\Theta_{N_1}-N'_1 = k-m, N'_{2,1} = \ell, N'_1=m] \\
	&= \IP[\Theta_{N_1}-N'_1 = k-m, \tilde{Y}_{k-m} = \ell, N'_1=m]\\
	&= \IP[\Theta_{N_1}-N'_1 = k-m, N'_1=m] \IP[\tilde{Y}_{k-m} = \ell]\\
	&= \IP[N'_{2,1} = \ell]\IP[\Theta_{N_1}= k]\IP[N'_1=m].
	\end{align*}
\end{proof}

Consequently, we obtain a coupling between~$\II^u$ and~$\MM^u$ (and therefore between $\I^u_{\sK}$ and $\M^u_{\sK}$). 
We will denote by~$\Upsilon$ the coupling event associated to this coupling (that is, $\Upsilon = \{\II^u=\MM^u\}$). In Section~\ref{Proofmaintheo}, we will obtain an upper bound for $\IP[\Upsilon^c]$.

\section{Toolbox}\label{toolbox}

In this section, we will prove some auxiliary results that will be needed to prove Theorem \ref{Maintheo}.
\medskip

Recall that we assume $\|\hat{x}\| ~\geq 4\diam(\sK_1)+3$, and we take~$R=\frac{\|\hat{x}\|-1}{2}$, so that $R\geq 2\diam(\sK_1)+1$. Furthermore, let~$\delta=\frac{\diam\displaystyle{(\sK_1)\vee 1}}{R}$, and recall that $q=\displaystyle\inf_{y\in \partial_e \sV_R} P_y[\tau_{\sK}=\infty]$.

{\lem 
\label{lemma_escape}
There exist $c_1>0$ and $0<\delta_1\leq \frac{1}{2}$, depending only on the dimension~$d$, such that, for~$\delta\leq \delta_1$, we have 
\begin{align*}
q \geq 1 - c_1 \frac{\capa(\emph{\sK}_1)}{R^{d-2}} \geq \frac{1}{2}.
\end{align*}
}

\begin{proof}
Consider~$y\in \partial_e \sV_R$. By Proposition 6.5.1 of~\cite{LL10}, since $R\geq 2\diam(\sK_1)$, we obtain that 
\begin{align}
\label{prop651}
P_y[\tau_{\sK_1}<\infty] \asymp \frac{\capa(\sK_1)}{\|y\|^{d-2}}.
\end{align}
Now observe that
\begin{align}
\label{Est}
P_y[\tau_{\sK}=\infty] \geq 1- P_y[\tau_{\sK_1}<\infty]-P_y[\tau_{\sK_2}<\infty].
\end{align}
From \eqref{prop651}, since $\|y\| ~\geq R$, we deduce that there exists a constant $\gamma_1>0$ depending only on $d$ such that
\begin{align*}
P_y[\tau_{\sK_1}<\infty] \leq \gamma_1\frac{\capa(\sK_1)}{R^{d-2}}.
\end{align*}
The same bound can be obtained for the term $P_y[\tau_{\sK_2}<\infty]$.
Plugging these two bounds into (\ref{Est}), we obtain that
\begin{align*}
\inf_{y\in \partial_e \sV_R}P_y[\tau_{\sK}=\infty] \geq 1- 2\gamma_1\frac{\capa(\sK_1)}{R^{d-2}}.
\end{align*}

Now, since~$\sK_1\subset \BB_{(2\diam(\sK_1))\vee 1}(0)$ we have $\capa(\sK_1)\leq \capa(\BB_{(2\diam(\sK_1))\vee 1}(0))$, and by Proposition 6.5.2 of \cite{LL10}, $\capa(\BB_{(2\diam(\sK_1))\vee 1}(0))\leq \gamma_2 (\diam(\sK_1)\vee 1)^{d-2}$, where $\gamma_2>0$ depends only on $d$. Finally, we deduce that there exists $0<\delta_1\leq \frac{1}{2}$ depending only on $d$ such that 
\begin{align*}
\inf_{y\in \partial_e \sV_R}P_y[\tau_{\sK}=\infty] \geq 1- 2 \gamma_1\gamma_2\delta^{d-2}\geq \frac{1}{2}
\end{align*}
for $\delta\leq \delta_1$. This concludes the proof of the lemma.
\end{proof}

{\lem 
\label{lemma1}  Consider $\delta_1$ from Lemma \ref{lemma_escape}. Suppose that $\delta\leq \delta_1$. In this case we have, for any~$y\in  \partial\emph{\sK}$, 
\begin{align*}
\bar{e}_{\emph{\sK}}(y) \geq \frac{1}{4} \bar{e}_{\emph{\sK}_1}(y).
\end{align*}
The result is also true when~$\bar{e}_{\emph{\sK}_1}(y)$ is replaced by~$\bar{e}_{\emph{\sK}_2}(y)$.
}

\begin{proof}
We prove the first assertion only, that is, the one involving~$\bar{e}_{\sK_1}(y)$. The case with~$\bar{e}_{\sK_2}(y)$ is analogous.

If~$y\in\partial \sK_2$ then~$\bar{e}_{\sK_1}(y)=0$ and the result follows trivially. Hence, let us suppose that~$y\in\partial \sK_1$.
Recall that we are denoting the union~$\sK_1\cup \sK_2$ simply by~$\sK$, and from the definition of the harmonic measure, 
\begin{align*}
\bar{e}_{\sK}(y) = \frac{P_y[\tau_{\sK}=\infty]}{\capa(\sK)} \geq \frac{P_y[\tau_{\sK}=\infty]}{2\capa(\sK_1)},
\end{align*}
where the inequality follows from the subadditivity property of the capacity.

Next, when dealing with the numerator, we obtain that
\begin{align*}
P_y[\tau_{\sK}=\infty] &= P_y[\tau_{\sK}=\infty, \tau_{\partial_e \BB_R^1}<\tau_{\sK_1}] \\
&= \sum_{x\in\partial_e \BB_R^1} P_y\Big[\tau_{\sK}=\infty, \tau_{\partial_e \BB_R^1}<\tau_{\sK_1}, X_{\tau_{\partial_e \BB_R^1}} = x\Big] \\
&= \sum_{x\in\partial_e \BB_R^1} P_x[\tau_{\sK}=\infty] P_y[\tau_{\partial_e \BB_R^1}<\tau_{\sK_1}, X_{\tau_{\partial_e \BB_R^1}} = x] \\
&\geq P_y[\tau_{\partial_e \BB_R^1}<\tau_{\sK_1}] \inf_{x\in\partial_e \BB_R^1} P_x[\tau_{\sK}=\infty],
\end{align*}
where we used the Markov property in the third equality. 
But, observe that the escape probability (or equilibrium measure) of~$\sK_1$ satisfies
\begin{align*}
e_{\sK_1}(y) \leq P_y[\tau_{\partial_e \BB_R^1}<\tau_{\sK_1}]
\end{align*}
for any~$y\in\partial \sK_1$, and also, by Lemma~\ref{lemma_escape},
\begin{align*}
\inf_{x\in\partial_e \BB_R^1} P_x[\tau_{\sK}=\infty]\geq \frac{1}{2}.
\end{align*}
Gathering these facts, we conclude the proof of the lemma.
\end{proof}

Recalling the transition density~$g$ defined in~\eqref{trans_g}, we prove the following

{\prop \label{Prop1} Consider $\delta_1$ from Lemma \ref{lemma_escape}. There exists a positive constant $c_2$ depending only on the dimension $d$, such that, for all $\delta\leq \delta_1$, we have that
\begin{align*}
|g(z',z) - 1| \leq c_2,
\end{align*}
for all~$z,z'\in\Sigma$.
\label{Prop_harmonic_aprox}
}

\begin{proof}
To simplify the notation in the proof, we denote by~$\rho_x$ the conditional probability of hitting~$\sK_1$ before~$\sK_2$, starting at~$x\in \partial_e \sV_R$ and given that the union of these sets is visited at a finite time, namely
\begin{align*}
\rho_x = P_{x}[\tau_{\sK_1}<\tau_{\sK_2}\mid \tau_{\sK}<\infty].
\end{align*}
Recall that, by definition, 
\begin{align*}
g(z',z) = \frac{P_{S_f(z')}[X_{\tau_{\sK}}=S_0(z)\mid \tau_{\sK}<\infty]}{\bar{e}_{\sK}(S_0(z))},
\end{align*}
and, without loss of generality, let us denote~$S_f(z')=x$ and~$S_0(z)=y$. We know that~$x\in \partial_e \sV_R$,
 and additionally let us suppose
 that~$y\in\partial \sK_1$. The case with~$y\in\partial \sK_2$ 
is analogous.

Then, we can state that
\begin{align}
\label{TY1}
P_{x}[X_{\tau_{\sK}}=y\mid \tau_{\sK}<\infty] &= P_{x}[X_{\tau_{\sK}}=y, \tau_{\sK_1}<\tau_{\sK_2}|\tau_{\sK}<\infty] \nonumber\\
&= \rho_x P_{x}[X_{\tau_{\sK}}=y\mid \tau_{\sK_1}<\tau_{\sK_2}, \tau_{\sK}<\infty].
\end{align}

But observe that the conditional probability multiplying the term~$\rho_x$ above is equal to
\begin{align}
\label{TY2}
P_{x}[X_{\tau_{\sK}}=y\mid\tau_{\sK_1}<\tau_{\sK_2}, \tau_{\sK_1}<\infty] = 
\frac{P_{x}[X_{\tau_{\sK}}=y\mid\tau_{\sK_1}<\infty]}{P_{x}[\tau_{\sK_1}<\tau_{\sK_2} \mid \tau_{\sK_1}<\infty]}.
\end{align}
Now, the denominator in the last expression can be written as
\begin{align}
\label{TY3}
\frac{P_{x}[\tau_{\sK_1}<\tau_{\sK_2} , \tau_{\sK_1}<\infty]}{P_{x}[\tau_{\sK_1}<\infty]} 
&=\frac{P_{x}[\tau_{\sK_1}<\tau_{\sK_2} , \tau_{\sK}<\infty]}{P_{x}[\tau_{\sK_1}<\infty]}\nonumber  \\
&= \frac{P_{x}[\tau_{\sK}<\infty]}{P_{x}[\tau_{\sK_1}<\infty]}  \rho_x.
\end{align}
Gathering (\ref{TY1}), (\ref{TY2})  and (\ref{TY3}), we obtain that 
\begin{align*}
P_{x}[X_{\tau_{\sK}}=y\mid \tau_{\sK}<\infty]
&= \frac{P_{x}[\tau_{\sK_1}<\infty]}{P_{x}[\tau_{\sK}<\infty]} P_{x}[X_{\tau_{\sK}}=y|\tau_{\sK_1}<\infty]\\
&\leq P_{x}[X_{\tau_{\sK}}=y|\tau_{\sK_1}<\infty\Big].
\end{align*}
Then,  using the fact that $y\in \partial \sK_1$ and applying Proposition 6.5.4 of~\cite{LL10} we obtain that, for $\delta\leq \delta_1$, 
\begin{align*}
P_{x}[X_{\tau_{\sK}}=y|\tau_{\sK_1}<\infty] \leq P_{x}\Big[X_{\tau_{\sK_1}}=y|\tau_{\sK_1}<\infty\Big] =  \bar{e}_{\sK_1}(y)\Big(1+O(\delta)\Big).
\end{align*}
Using this last inequality and Lemma \ref{lemma1}, we obtain
\begin{align*}
P_{x}[X_{\tau_{\sK}}=y|\tau_{\sK}<\infty]
&\leq 4(1+O(\delta)) \bar{e}_{\sK}(y).
\end{align*}
From this last inequality, we easily deduce the proposition.
\end{proof}

Recalling~$N_1$ and~$N_{2,2}$ from Section~\ref{coupling}, we prove the following

{\lem 
	\label{N1N2null}
	Consider~$\delta_1$ from Lemma \ref{lemma_escape}. There exists a positive constant~$c_3$ depending only on the dimension $d$, such that, for all~$\delta\leq \delta_1$, we have that 
	\begin{align*}
	\IP[N_1\geq 1, N_{2,2} =0] \leq c_3 \sqrt{u} \frac{\capa(\emph{\sK}_1)^{\frac{3}{2}}}{R^{d-2}}.
	\end{align*}	
}

\begin{proof}
Let us denote~$\theta=u\capa(\sK)$. Recall that~$N_1$ and~$N_{2,2}$ are independent Poisson random variables, with parameters~$(1-q)\theta$ and~$\frac{q}{2}\theta$, respectively. Thus 
\begin{align*}
\IP[N_1 \geq 1, N_{2,2} =0] = (1-e^{-(1-q)\theta}) e^{-\frac{q}{2}\theta} \leq \frac{(1-q)\theta}{\sqrt{\frac{q}{2}\theta}} = \frac{\sqrt{2}(1-q)}{\sqrt{q}}\sqrt{\theta}, 
\end{align*}
where we used the facts that~$1-e^{-x}\leq x$ for~$x\geq 0$, and~$e^{-x}\leq x^{-1/2}$ for~$x>0$. The proof is concluded by using the fact that $\capa(\sK)\leq 2\capa(\sK_1)$, and then applying Lemma~\ref{lemma_escape} to the right-hand term of the above display.
\end{proof}

{\lem 
	\label{N2_exp}
	Consider $\delta_1$ from Lemma \ref{lemma_escape}. There exists a positive universal constant~$c_4$ such that, for all~$\delta\leq \delta_1$, we have that
	\begin{align*}
	\IE[N_{2,2}^{-1/2}\ind_{\{N_{2,2}\geq 1\}}] \leq c_4 (u\capa(\emph{\sK}_1))^{-\frac{1}{2}}.	
	\end{align*}	
}

\begin{proof}
 Again, let us denote~$\theta=u\capa(\sK)$. Using Cauchy-Schwarz inequality and the fact that~$6x^2\geq (x+1)(x+2)$ for~$x\geq 1$, we have that
 \begin{align*}
 \IE[N_{2,2}^{-1/2}\ind_{\{N_{2,2}\geq 1\}}] &\leq \sqrt{\IE[N_{2,2}]\IE[N_{2,2}^{-2}\ind_{\{N_{2,2}\geq 1\}}]} \\
 &\leq \sqrt{\frac{q}{2}\theta} \sqrt{6\IE\Big[\frac{1}{(N_{2,2}+2)(N_{2,2}+1)}\ind_{\{N_{2,2}\geq 1\}}\Big]} \\
 &\leq \sqrt{3} \sqrt{\theta} \sqrt{\IE\Big[\frac{1}{(N_{2,2}+2)(N_{2,2}+1)}\Big]},
 \end{align*}
 and since~$\IE[1/[(N_{2,2}+2)(N_{2,2}+1)]] \leq (\frac{q}{2}\theta)^{-2}$, we obtain the result with~$c_4=4\sqrt{3}$ after using the bound~$q\geq 1/2$ from Lemma~\ref{lemma_escape}, and the fact that $\capa(\sK)\geq \capa(\sK_1)$.
 \end{proof}

{\lem 
\label{Lem_Poi_TV}  
If $X$ is a Poisson distributed random variable with parameter~$\theta$, and $X_k=X+k$ for~$k\in\Z$, then
\begin{align*}
\dtv( X , X_k ) \leq \frac{|k|}{\sqrt{\theta}}.
\end{align*}
}

\begin{proof}
Let~$(Y_k)_{k\in\Z}$ be a family of Poisson($\theta$) distributed random variables and define~$Z_k=Y_k+k$, for~$k\in\Z$. It is elementary to see that~$\dtv( Z_i , Z_{i+1} ) = \dtv( X , Z_1 )$, for any~$i\in\Z$.
Now, for any positive integer~$k$, 
\begin{align*}
\dtv( X , Z_k ) \leq \dtv( X , Z_1 ) + \sum_{j=1}^{k-1}\dtv( Z_j , Z_{j+1} ) = k\dtv( X , Z_1 ),
\end{align*}
and for any negative integer~$k$,
\begin{align*}
\dtv( X , Z_k ) \leq \dtv( X , Z_{-1} ) + \sum_{j=1}^{|k|-1}\dtv( Z_{-j} , Z_{-(j+1)} ) = |k|\dtv( X , Z_1).
\end{align*}
But, from Lemma 1 of~\cite{BP15}, one has 
\begin{align*}
\dtv( X , Z_1 ) \leq \frac{1}{\sqrt{\theta}},
\end{align*}
and this concludes the proof, since the result is obvious for~$k=0$.
\end{proof}

\section{Proof of Theorem \ref{Maintheo}}
\label{Proofmaintheo}
In this section we will show that the coupling we constructed 
in Section~\ref{coupling} between $\II^u$ and $\MM^u$ 
is successful with high probability, when $\|\hat{x}\|$ is large. Recall that we denoted
 by~$\Upsilon$ the coupling event of~$\II^u$ and~$\MM^u$, and by~$\IP$ the probability
on the space they are jointly constructed. 
This will automatically give an upper bound on the total 
variation distance between the two processes $\I^u_{\sK}$ and $\M^u_{\sK}$, since 
\begin{align*}
\dtv(\I^u_{\sK},\M^u_{\sK}) \leq \dtv(\II^u,\MM^u) \leq \IP[\Upsilon^c].
\end{align*}

The goal of this section is to estimate $\IP[\Upsilon^c]$ from above. 
Before starting, we mention that we will use the notation from 
Section~\ref{coupling}.


As a first step to obtain an upper bound for $\IP[\Upsilon^c]$, we consider the events
\begin{align}
\label{events_A}
\sA_i^{k,\ell} = \Bigg\{1\wedge\sup_{z\in\Sigma} |\Psi(z)-1| \leq (1+i)\frac{k}{\ell}\Bigg\}, 	
\end{align}
for~$i,k,\ell\in\IN$ with~$\ell\geq 4$, and prove the following

{\prop 
\label{Prop_Ac}
Consider $\delta_1$ from Lemma \ref{lemma_escape}. There exists a positive constant~$c_5$ depending only on the dimension $d$, such that, for all $\delta\leq \delta_1$ and for~$i,k,\ell\in\IN$, $\ell\geq 4$, we have that
\begin{align*}
\IP\Big[(\emph{\sA}_i^{k,\ell})^c \mid N_1=k, N_{2,2}=\ell\Big] \leq \frac{c_5}{(1+i)^3}.
\end{align*}
}

\begin{proof}
	

Let~$i,k,\ell\in\IN$, $\ell\geq 4$. Also, consider the event 
\begin{align*}
\sH := \Bigg\{\sup_{z,z'\in\Sigma} |G^{\II}_{\Theta_{N_1}}(z)-G^{\II}_{\Theta_{N_1}}(z')| \leq \frac{\Xi_{2,2}}{2}\Bigg\}.
\end{align*}
First, observe that
\begin{align*}
\IE\Big[1\wedge&\sup_{z\in\Sigma} |\Psi(z)-1|^3 ~\Big|~ N_1=k, N_{2,2}=\ell \Big] \\
&\leq \IE\Big[1\wedge\sup_{z\in\Sigma} |\Psi(z)-1|^3\ind_{\sH} ~\Big|~ N_1=k, N_{2,2}=\ell \Big] + \IP[\sH^c ~|~  N_1=k, N_{2,2}=\ell].
\end{align*}
Then, for $\delta\leq \delta_1$, we bound the first term of the last sum from above by
\begin{align*}
2^3 \IE\Bigg[&\Bigg( \frac{\sup_{z,z'\in\Sigma} |G^{\II}_{\Theta_{N_1}}(z)-G^{\II}_{\Theta_{N_1}}(z')|}{\Xi_{2,2}} \Bigg)^3 ~\Big|~ N_1=k, N_{2,2}=\ell \Bigg] \\
&\leq 2^6\IE\Bigg[\sup_{z\in\Sigma} \Bigg|\sum_{j=1}^{k} \sum_{n=1}^{T_j} \xi_{\Theta_{j-1}+n}  \Big(1 - \tilde{g}(z_{\Theta_{j-1}+n-1},z)\Big)\Bigg|^3 \Bigg(\sum_{i=1}^{\ell}\xi_{\Theta_k+N_{2,1}+i}\Bigg)^{-3} \Bigg] \\
&\leq 2^6 c_2^3\IE\Bigg[\Bigg(\sum_{j=1}^{k} \sum_{n=1}^{T_j} \xi_{\Theta_{j-1}+n} \Bigg)^3\Bigg] \IE\Bigg[\Bigg(\sum_{i=1}^{\ell}\xi_{\Theta_k+N_{2,1}+i}\Bigg)^{-3} \Bigg],
\end{align*}
where in the first inequality we used the fact that 
\begin{align*}
\sup_{z,z'\in\Sigma} |G^{\II}_{\Theta_{N_1}}(z)-G^{\II}_{\Theta_{N_1}}(z')| ~\text{ and }~ \Xi_{2,2}
\end{align*}
are conditionally independent given $N_1$ and $N_{2,2}$, and in the last inequality we used Proposition~\ref{Prop1}.
Next, observe that
\begin{align*}
\IE\Bigg[\Bigg(\sum_{j=1}^{k} \sum_{n=1}^{T_j} \xi_{\Theta_{j-1}+n} \Bigg)^3\Bigg] \leq \Bigg(\sum_{j=1}^{k} \Big\|\sum_{n=1}^{T_j} \xi_{\Theta_{j-1}+n}\Big\|_3 \Bigg)^3 = k^3 \Big\|\sum_{n=1}^{T_1} \xi_n\Big\|_3^3 \leq  k^3\gamma_1, 
\end{align*}
where~$\gamma_1$ is a positive constant and the last inequality can be obtained in an elementary way by observing that~$T_1\preceq 1+\text{Geo}(q)$. On the other hand, since~$(\sum_{i=1}^{\ell}\xi_{\Theta_k+N_{2,1}+i})^{-1}$ has Inverse Gamma distribution with parameters~$(\ell, 1)$, we have
\begin{align*}
\IE\Bigg[\Bigg(\sum_{i=1}^{\ell}\xi_{\Theta_k+N_{2,1}+i}\Bigg)^{-3} \Bigg]= \frac{1}{(\ell-1)(\ell-2)(\ell-3)} \leq \frac{\gamma_2}{\ell^{3}},
\end{align*}
for a positive constant~$\gamma_2$.
Thus, it holds that
\begin{align*}
\IE\Big[1\wedge\sup_{z\in\Sigma} |\Psi(z)-1|^3\ind_{\sH} ~\Big|~ N_1=k, N_{2,2}=\ell \Big] \leq 2^6 c_2^3\gamma_1\gamma_2 \frac{k^3}{\ell^3}.
\end{align*}

Using the above estimates and Markov's inequality we obtain 
\begin{align*}
\IP[\sH^c ~|~  N_1=k, N_{2,2}=\ell] \leq 2^6 c_2^3\gamma_1\gamma_2 \frac{k^3}{\ell^3}.
\end{align*}
Finally, we obtain
\begin{align*}
\IE\Big[1\wedge&\sup_{z\in\Sigma} |\Psi(z)-1|^3 ~\Big|~ N_1=k, N_{2,2}=\ell \Big] \leq \gamma_3\frac{k^3}{\ell^3},
\end{align*}
and applying Markov's inequality once again, we have
\begin{align*}
\IP\Big[(\emph{\sA}_i^{k,\ell})^c \mid N_1=k, N_{2,2}=\ell\Big] \leq \frac{c_5}{(1+i)^3},
\end{align*}
for a positive constant $c_5$.
\end{proof}


Now, for~$k,\ell\in\IN$ with~$\ell\geq 4$, we define~$\sC_1^{k,\ell}=\sA_1^{k,\ell}$, $\sC_{i+1}^{k,\ell}=\sA_{i+1}^{k,\ell}\setminus\sA_i^{k,\ell}$, for~$i\geq 1$.
Also, let us introduce the event
\begin{align}
\label{eventD}
\sD := \{\Theta_{N_1}+N_{2,1} = N'_1 + N'_{2,1}\}. 
\end{align}
The event $\sD$ corresponds to the coupling event of the cardinalities of $\II^u$ and $\MM^u$. 

We start by decomposing $\IP[\Upsilon^c]$ in the following way:
\begin{align}
\label{decomposit}
\IP[\Upsilon^c]= &~ \IP[\Upsilon^c\cap\sD^c] + \IP[\Upsilon^c\cap\sD\cap\{N_1=0\}] \nonumber \\
 &+\IP[\Upsilon^c\cap\sD\cap\{ N_1\geq 1, N_{2,2}=0\}] \nonumber\\
 &+ \sum_{k\geq 1}\sum_{\ell\geq 1} \IP[\Upsilon^c\cap\sD\cap \{N_1=k, N_{2,2}=\ell\}].
\end{align}
As we will see below, the main contributions to obtain the exponent $3/2$ for $\capa(\sK_1)$ in Theorem~\ref{Maintheo} will be given by the first and the last terms in \eqref{decomposit}.

First, observe that, by the coupling construction from Section~\ref{coupling}, we automatically have
\begin{align}
\IP[\Upsilon^c\cap\sD\cap\{N_1=0\}] = 0.
\label{null_prob}
\end{align}

Next, since~$\sD\cap\sC_i^{k,\ell}\cap \{N_1=k\} \cap \{N_{2,2}=\ell\}$, for~$i,k,\ell\in\IN$, $\ell\geq 4$, are $\sigma(W)$-measurable (recall that $W$ was introduced in Section \ref{coupling}), we have
\begin{align*}
 \IP[\Upsilon^c\cap\sD\cap\sC_i^{k,\ell} \cap\{N_1=k, N_{2,2}=\ell\}] = \IE\Big[{\bf 1}_{\sD\cap\sC_i^{k,\ell}\cap \{N_1=k\} \cap \{N_{2,2}=\ell\}}\IP[\Upsilon^c\mid W]\Big],
\end{align*}
for~$i,k,\ell\in\IN$, $\ell\geq 4$. Then, again by the coupling construction from Section~\ref{coupling}, for~$k,\ell\in\IN$, $\ell\geq 4$ we have, on~$\sD\cap \{N_1=k\} \cap \{N_{2,2}=\ell\}$,
\begin{align*}
\IP[\Upsilon^c\mid W] \leq \dtv\Big(\IP[{\bf Y}\in\cdot\mid W] , \IP[{\bf Y'}\in\cdot\mid W] \Big).
\end{align*}
Therefore, using Proposition 5.1 of~\cite{BGP16} (choosing~$\delta_0=1$ in that proposition), for~$k,\ell\in\IN$, $\ell\geq 4$, we obtain, on the sets~$\sD\cap\sC_i^{k,\ell}\cap \{N_1=k\} \cap \{N_{2,2}=\ell\}$,
\begin{align*}
\IP[\Upsilon^c\mid W] \leq \Bigg(\gamma_1(1+i)\frac{k}{\sqrt{\ell}}\Bigg)\wedge 1 \leq \gamma_1(1+i)\frac{k}{\sqrt{\ell}},
\end{align*}
where~$\gamma_1$ is a constant greater than $1$. Hence we obtain, for $\delta\leq \delta_1$,
\begin{align*}
\sum_{k\geq 1}&\sum_{\ell\geq 4} \IP[\Upsilon^c\cap\sD\cap\{N_1=k, N_{2,2}=\ell\}] \\
&=\sum_{k\geq 1}\sum_{\ell\geq 4}\sum_{i\geq 1} \IP[\Upsilon^c\cap\sD\cap\sC_i^{k,\ell} \cap\{N_1=k, N_{2,2}=\ell\}] \\ 
&\leq \sum_{k\geq 1}\sum_{\ell\geq 4}\sum_{i\geq 1} \gamma_1(1+i)\frac{k}{\sqrt{\ell}}\IP[\sC_{i}^{k,\ell}\mid N_1=k, N_{2,2}=\ell]\IP[N_1=k]\IP[ N_{2,2}=\ell] \\
&\leq \gamma_1\sum_{k\geq 1}\sum_{\ell\geq 4}\Bigg[\frac{k}{\sqrt{\ell}}\IP[N_1=k]\IP[ N_{2,2}=\ell] \\
&\hspace{4cm}  \times\Big(2+\sum_{i\geq 2} (1+i) \IP[(\sA_{i-1}^{k,\ell})^c\mid N_1=k, N_{2,2}=\ell]\Big)\Bigg] \\
&\leq \gamma_2 \sum_{k\geq 1}\sum_{\ell\geq 4}\frac{k}{\sqrt{\ell}}\IP[N_1=k]\IP[ N_{2,2}=\ell],
\end{align*}
where we used Proposition~\ref{Prop_Ac} in the last inequality.

Now, using the last bound and observing that 
\begin{align*}
\IP[\Upsilon^c\cap\sD\cap\{N_1\geq 1,& 1\leq N_{2,2}\leq 3\}] \\
&=\sum_{k\geq 1} \sum_{\ell=1}^{3}\IP[\Upsilon^c\cap\sD\cap\{N_1=k, N_{2,2}=\ell\}] \\
&\leq 2\sum_{k\geq 1} \sum_{\ell=1}^{3} \frac{k}{\sqrt{\ell}}\IP[N_1=k]\IP[ N_{2,2}=\ell],
\end{align*}
we deduce that, for $\delta\leq \delta_1$,
\begin{align}
\label{lbound}
\sum_{k\geq 1}\sum_{\ell\geq 1} \IP[\Upsilon^c\cap\sD\cap& \{N_1=k, N_{2,2}=\ell\}] \nonumber\\
&=\sum_{k\geq 1}\sum_{\ell\geq 4}\sum_{i\geq 1} \IP[\Upsilon^c\cap\sD\cap\sC_i^{k,\ell} \cap\{N_1=k, N_{2,2}=\ell\}]  \nonumber\\
&\hspace{0.5cm}+ \IP[\Upsilon^c\cap\sD\cap\{N_1\geq 1, 1\leq N_{2,2}\leq 3\}] \nonumber\\
& \leq \gamma_3 \sum_{k\geq 1}\sum_{\ell\geq 1} \frac{k}{\sqrt{\ell}}\IP[N_1=k]\IP[ N_{2,2}=\ell] \nonumber\\
&= \gamma_3 \IE[N_1]\IE[N_{2,2}^{-1/2}\ind_{\{N_{2,2}\geq 1\}}] \nonumber\\
&\leq \gamma_4 \sqrt{u} \frac{\capa(\sK_1)^{\frac{3}{2}}}{R^{d-2}},
\end{align}
where we used Lemmas~\ref{lemma_escape} and~\ref{N2_exp} to obtain the last inequality.

We just bound the term $\IP[\Upsilon^c\cap\sD^c]$ in \eqref{decomposit} by $\IP[\sD^c]$, and then we use the following

{\lem 
\label{lemma_D}
Consider~$\delta_1$ from Lemma \ref{lemma_escape}. There exists a positive constant~$c_6$ depending only on the dimension $d$, such that, for all~$\delta\leq \delta_1$, we have that
\begin{align*}
\IP[\emph{\sD}^c] \leq c_6 \sqrt{u} \frac{\capa(\emph{\sK}_1)^{\frac{3}{2}}}{R^{d-2}}.	
\end{align*}	
}

\begin{proof}
Recall that ~$\sD^c = \{\Theta_{N_1}+N_{2,1} \neq N'_1 + N'_{2,1}\}$.
From the definition of~$N_{2,1}$ and~$N'_{2,1}$ given in~\eqref{def_N21} we have that, for any~$k\in\Z$,
\begin{align*}
\IP[\Theta_{N_1}+N_{2,1} \neq N'_1 + N'_{2,1} \mid & \Theta_{N_1} - N'_1 = k] \\
&= \IP[N'_{2,1} \neq k + N_{2,1} \mid \Theta_{N_1} - N'_1 = k]\\
&= \IP[\tilde{Y}_k \neq k + \tilde{X}_k] ,
\end{align*}
which implies that
\begin{align*}
\IP[\Theta_{N_1}+N_{2,1} \neq N'_1 + N'_{2,1}] &= \sum_{k\in\Z} \IP[\tilde{Y}_k \neq k + \tilde{X}_k] \IP[\Theta_{N_1} - N'_1 = k] \\
& = \sum_{k\in\Z} \dtv\Big( \tilde{Y}_k , (k + \tilde{X}_k) \Big) \IP[\Theta_{N_1} - N'_1 = k].
\end{align*}
But, from Lemma~\ref{Lem_Poi_TV},
\begin{align*}
\dtv\Big( \tilde{Y}_k , (k + \tilde{X}_k) \Big) \leq |k| \frac{\sqrt{2}}{\sqrt{q}}  (u\capa(\sK))^{-\frac{1}{2}} \leq 
2|k| (u\capa(\sK))^{-\frac{1}{2}},
\end{align*}
where we used the fact that~$q\geq 1/2$. Thus, since $\IE[\Theta_{N_1}]=\IE[N'_1]$, we have that
\begin{align*}
\IP[\Theta_{N_1}+N_{2,1} \neq N'_1 + N'_{2,1}] &\leq 2(u\capa(\sK))^{-\frac{1}{2}} \IE\Big[|\Theta_{N_1} - N'_1|\Big]\\
& \leq 4(u\capa(\sK))^{-\frac{1}{2}} \IE[\Theta_{N_1}].
\end{align*}	
Recalling that $\Theta_{N_1}=\sum_{j=1}^{N_1}T_j$, the fact that the random variables $(T_j)_{j\geq 1}$ are independent of $N_1$ and that $T_j\preceq 1+\text{Geo}(q)$ for all $j\geq 1$, we obtain
\begin{align*}
\IE[\Theta_{N_1}] \leq 3(1-q)u\capa(\sK).
\end{align*}
Finally, the result follows by just applying Lemma~\ref{lemma_escape} to bound~$1-q$ from above by $c_1\frac{\capa\big(\sK_1\big)}{R^{d-2}}$, and using that $\capa(\sK) \leq 2\capa(\sK_1)$.
\end{proof}

Thus, gathering~\eqref{decomposit}, \eqref{null_prob}, \eqref{lbound} and Lemmas~\ref{N1N2null} and~\ref{lemma_D}, we conclude that there exists a positive constant~$\gamma_5$, depending only on the dimension, such that, for $\delta \leq \delta_1$,
\begin{align*}
\IP[\Upsilon^c]\leq \gamma_5 \sqrt{u} \frac{\capa(\sK_1)^{\frac{3}{2}}}{R^{d-2}}.
\end{align*}

Finally, using the fact that $\dist(\sK_1,\sK_2)\leq 3R$, we deduce that there exist positive constants~$\gamma_6$ and $\gamma_7$ depending only on the dimension such that, for $\dist(\sK_1,\sK_2)\geq \gamma_6 (\diam(\sK_1)\vee 1)$,
\begin{align*}
\IP[\Upsilon^c]\leq \gamma_7 \sqrt{u} \frac{\capa(\sK_1)^{\frac{3}{2}}}{\dist(\sK_1,\sK_2)^{d-2}}.
\end{align*}
Observe that, in the case the set $\sK_1$ is a single point in~$\Z^d$, using~\cite{DRS14}, Claim 2.5, the condition $\dist(\sK_1,\sK_2)\geq \gamma_6 (\diam(\sK_1)\vee 1)$ can be relaxed to $\dist(\sK_1,\sK_2)\geq \gamma_6 \diam(\sK_1)=0$. This concludes the proof of Theorem~\ref{Maintheo}.
\qed

\section{Proof of Corollary \ref{Coro1}} 
\label{ProofCoro1}
Corollary~\ref{Coro1} is a direct consequence of Theorem~\ref{Maintheo} and the next proposition.

{\prop \label{Propo1} Suppose that we are given two functions $f_1:\{0,1\}^{\emph{\sK}_1}\to [0,1]$ and $f_2:\{0,1\}^{ \emph{\sK}_2}\to [0,1]$ that depend only on the configuration of the random interlacements inside the sets $\emph{\sK}_1$ and $\emph{\sK}_2$, respectively. We have that
\begin{align*}
|\emph{Cov}(f_1(\I^u_{\emph{\sK}_1}), f_2(\I^u_{\emph{\sK}_2}))| \leq 3 \dtv(\I^u_{\emph{\sK}},\M^u_{\emph{\sK}}).
\end{align*}
}

\begin{proof}
We first introduce $\M^u_{\sK_i}:=\M^u_{\sK}\cap \sK_i$, for $i=1,2$. By definition of the covariance, we have that
\begin{equation*}
|\text{Cov}(f_1(\I^u_{\sK_1}), f_2(\I^u_{\sK_2}))|=|E[f_1(\I^u_{\sK_1})f_2(\I^u_{\sK_2})]-E[f_1(\I^u_{\sK_1})]E[f_2(\I^u_{\sK_2})]|.
\end{equation*}
Now, let $\I'_{\sK_2}$ be a copy of $\I^u_{\sK_2}$,
such that $\I'_{\sK_2}$ is independent of $\I^u_{\sK_1}$. We obtain that
\begin{align*}
\label{Covfin}
|\text{Cov}(f_1(\I^u_{\sK_1}), f_2(\I^u_{\sK_2}))|
&=|E[f_1(\I^u_{\sK_1})f_2(\I^u_{\sK_2})]-E[f_1(\I^u_{\sK_1})f_2(\I'_{\sK_2})]|\nonumber\\
&\leq | E[f_1(\I^u_{\sK_1})f_2(\I^u_{\sK_2})]-E[f_1(\M^u_{\sK_1})f_2(\M^u_{\sK_2})]|\nonumber\\
&\phantom{**}+|E[f_1(\M^u_{\sK_1})f_2(\M^u_{\sK_2})]-E[f_1(\I^u_{\sK_1})f_2(\I'_{\sK_2})]|.
\end{align*}
Next, since $f_1$ and $f_2$ are $[0,1]$-valued, it is not difficult to see that
\begin{align*}
| E[f_1(\I^u_{\sK_1})f_2(\I^u_{\sK_2})]-E[f_1(\M^u_{\sK_1})f_2(\M^u_{\sK_2})]| \leq \dtv(\I^u_{\sK},\M^u_{\sK}).
\end{align*}
 Also, note that
\begin{align*}
|E[f_1(\M^u_{\sK_1})f_2(\M^u_{\sK_2})]-E[f_1(\I^u_{\sK_1})f_2(\I'_{\sK_2})]|\leq \dtv(\M^u_{\sK},\I^u_{\sK_1}\cup\I'_{\sK_2}).
\end{align*}
Then we observe that, by definition of the NS process, $\M^u_{\sK_1}$ and $\M^u_{\sK_2}$ are independent. By symmetry, this leads to
\begin{align*}
 \dtv(\M^u_{\sK},\I^u_{\sK_1}\cup\I'_{\sK_2})\leq 2  \dtv(\M^u_{\sK_1},\I^u_{\sK_1}).
\end{align*}
Combining this last inequality with the fact that
\begin{align*}
\dtv(\M^u_{\sK_1},\I^u_{\sK_1})\leq \dtv(\I^u_{\sK},\M^u_{\sK}),
\end{align*}
we obtain Proposition~\ref{Propo1}.
\end{proof}

\section*{Acknowledgements}
\noindent
Diego F.~de Bernardini was partially supported by S\~ao Paulo Research Foundation (FAPESP) (grant 2014/14323-9). Christophe Gallesco was partially supported by FAPESP (grant 2017/19876-4) and CNPq (grant 312181/2017-5). Serguei Popov was partially supported by CNPq (grant 300886/2008--0). The three authors were partially supported by FAPESP (grant 2017/02022-2 and grant 2017/10555-0).

\end{document}